\documentclass[a4paper,11pt]{article}

\usepackage{amssymb}
\usepackage[T1]{fontenc}
\usepackage[utf8]{inputenc}
\usepackage{mathrsfs}
\usepackage{amsfonts,amsmath,lmodern,lastpage,graphicx,nccfoots,caption}
\usepackage{afterpage}
\usepackage{epsfig}
\usepackage{rlepsf}
\usepackage{graphicx}
\graphicspath{{./Figures/}}

\usepackage[all]{xy}
\usepackage{amscd}
\usepackage{comment}

\newtheorem{Theorem}{Theorem}[section]

\newtheorem{Definition}[Theorem]{Definition}

\newtheorem{Proposition}[Theorem]{Proposition}
\newtheorem{Corollary}[Theorem]{Corollary}
\newtheorem{Remark}[Theorem]{Remark}

\newenvironment{Proof}[1][Proof]{\textbf{#1.} }{\ \rule{0.5em}{0.5em}}

\def \tr {\triangleright}

\def \G {\mathcal{G}}

\usepackage{hyperref}

\begin{document}

\thispagestyle{plain}



\begin{center}

\Large

\textsc{Invariants and TQFT's for cut cellular surfaces from finite 2-groups}

\end{center}



\begin{center}

\textit{Diogo Bragan\c ca}  \smallskip \\

\begin{tabular}{l}
\small Multidisciplinary Centre for Astrophysics (CENTRA) \\
\small  Physics Department \\
\small Instituto Superior T\'ecnico, Universidade de Lisboa  \\
\small Av. Rovisco Pais, 1                    \\
\small 1049-001 Lisboa, Portugal       \\
\small e-mail: \texttt{diogo.braganca@tecnico.ulisboa.pt}   \\ 
\end{tabular}

\end{center}

\begin{center}

\textit{Roger Picken} \smallskip \\

\begin{tabular}{l}
\small Center for Mathematical Analysis, Geometry \\ 
\small and Dynamical Systems (CAMGSD), \\ 
\small Mathematics Department, \\
\small  Instituto Superior T\'ecnico, Universidade de Lisboa \\
\small Av. Rovisco Pais, 1 \\
\small 1049-001 Lisboa, Portugal \\
\small e-mail: \texttt{roger.picken@tecnico.ulisboa.pt}
\end{tabular}

\end{center}


\noindent
\textbf{Abstract} \\
In this brief sequel to a previous article, we recall the notion of a cut cellular surface (CCS), being a surface with boundary, which is cut in a specified way to be represented in the plane, and is composed of 0-, 1- and 2-cells. We obtain invariants of CCS's under Pachner-like moves on the cellular structure, by counting colourings of the 1- and 2-cells  with elements of a finite 2-group, subject to a ``fake flatness'' condition for each 2-cell. These invariants, which extend Yetter's invariants to this class of surfaces, are also described in a TQFT setting. A result from the previous article concerning the commuting fraction of a group is generalized to the 2-group context.

\medskip


\medskip



\noindent\textbf{keywords:} Cut cellular surface, TQFT, finite group, crossed module, 2-group, commuting fraction


\newpage

\section{Introduction}

In our previous work \cite{bp1} we studied invariants of a class of surfaces with boundary, obtained by counting certain $G$-colourings of the 1-cells of the surface, where $G$ is a finite group. We called these surfaces \emph{cut cellular surfaces} (CCS's), since they come equipped with a planar representation which arises from cutting the surface along some 1-cells to get a simply-connected planar region made up of 2-cells bounded by a circuit of 1- and 0-cells (see Section \ref{sec:ccs} for the full definition). Such surfaces include triangulated surfaces, but allow for a considerably more economical description in terms of the number of cells needed.  For instance, a triangulation of the 2-sphere $S$ requires at least four 2-cells, six 1-cells and four 0-cells, whereas its minimal representation as a CCS has just one 2-cell, one 1-cell and two 0-cells (see Section \ref{sec:ccs}). 

The invariants that we studied in \cite{bp1} involved counting the number of so-called flat $G$-colourings of the 1-cells of the surface, i.e. assignments of elements of $G$ to the 1-cells of the surface, such that, taking into account the orientation of the 1-cells, their product around the boundary of each 2-cell equals $1_G$, the identity element of $G$. For a triangulated surface without boundary, these invariants coincide with the Dijkgraaf-Witten invariants of the surface \cite{dw}. They are invariant under simple moves on the cellular structure, namely subdividing or combining 1-cells and subdividing or combining 2-cells. We showed that these two types of move generate the well-known Pachner moves for triangulated surfaces. The invariants also behave well under gluing of surfaces along shared boundary components and we showed that they give rise to a topological quantum field theory (TQFT). See the second section of our previous article \cite{bp1} for an introduction to the notion of TQFT.

The number of flat $G$-colourings for minimal CCS representations of some elementary surfaces, like the sphere, cylinder, pants surface and torus, has a group-theoretical significance, e.g. for the torus this is the number of commuting pairs of elements of $G$. Using topological arguments we were able to derive some group-theoretical properties, such as:
\begin{Proposition}
The number of conjugacy classes of $G$ is equal to the commuting fraction of $G$ times the order of $G$. 
\label{cf-prop}
\end{Proposition}
We recall that the commuting fraction of $G$ is defined to be the number of commuting pairs of elements of $G$ divided by the overall number of pairs.

The constructions in \cite{bp1} were intended to pave the way for an analogous approach using finite 2-groups, which is the subject of the present article. In Section \ref{sec:inv}, we recall the definition of a finite 2-group $\G$, also known as a finite crossed module. It consists of two finite groups $G$ and $H$, a group homomorphism from $H$ to $G$, and a left action of $G$ on $H$ by automorphisms, subject to two conditions.

We then define invariants of CCS's (Def. \ref{def:2grpinv}) which involve counting the number of $\G$-colourings of the surface, i.e. assignments of elements of $G$ to the 1-cells and elements of $H$ to the 2-cells. These assignments,  are subject to a ``fake flatness'' condition, which reduces to the flatness condition when the group $H$ is trivial. We prove several properties of the expressions of Definition \ref{def:2grpinv}, in particular that they are invariant under the aforementioned two types of move on the cellular structure. For triangulated surfaces these invariants correspond to Yetter's invariants \cite{yetter, FMPo}. In section \ref{sec:exp}  we calculate the invariant for some elementary examples. 

In Section \ref{sec:glue}, we describe how the invariant behaves when gluing two CCS's together along a common boundary component, and use this to get a TQFT for these surfaces. We focus on properties of the invariant for the cylinder, and in Proposition \ref{prop:gcf} we obtain a generalization of Proposition \ref{cf-prop}  in the 2-group context. Finally, in the conclusions of section \ref{sec:conc}, we comment on some features of the TQFT and give an interpretation for the invariants in terms of the notion of groupoid cardinality.

To make this article self-contained, we have repeated some material from \cite{bp1}. We invite the reader to consult this previous article for fuller details concerning a number of points.

To conclude this introduction we will say a brief word about notation. When we wish to describe a linear map $Z: V\rightarrow W$ in concrete terms, we may introduce a basis 
$\left\{e_i\right\}_{i=1,\dots,n}$ of $V$ and a basis $\left\{f_j\right\}_{j=1,\dots,m}$ of $W$. Then $Z$ is represented by an $m\times n$ matrix 
$[c_{ji}]$, where 
$$Z(e_i) =\sum_{j=1}^m c_{ji}f_j.
$$
We will be using the suggestive physicists' notation for the matrix elements $c_{ji}$, namely:
$$
c_{ji}=  \left \langle f_j\left | Z \right | e_i  \right \rangle.
$$

\vskip 0.3cm


\section{Cut cellular surfaces}

\label{sec:ccs}

We will be considering surfaces with boundary, which are cut in a specified way to be represented in the plane (like the well-known rectangle with opposite edges identified representing the torus), and which are composed of 0-, 1- and 2-cells, generalizing the familiar notion of a triangulated surface. 

\begin{figure}[htbp]
\centering
\includegraphics[width=10cm]{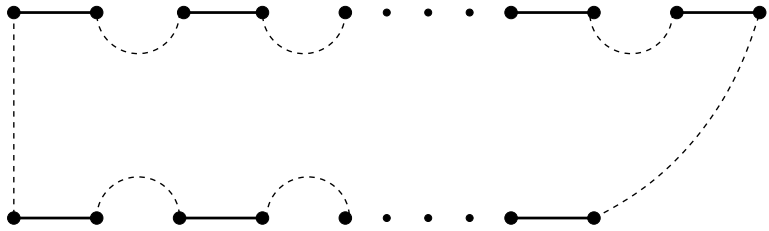}
\caption{General appearance of a cut cellular surface (CCS)}
\label{fig:ccs}
\end{figure}

\begin{Definition}

A cut cellular surface (CCS) is an orientable 2-manifold $M$ with boundary, endowed with a finite cell-structure, such that

\begin{itemize}

\item[a)] Each boundary component of $M$ consists of a single 0-cell and a single 1-cell.

\item[b)] $M$ has a specified planar representation, obtained by cutting $M$ along 1-cells in such a way as to obtain a simply connected region in the plane. The cut 1-cells are labeled and given an orientation to make explicit how they are identified in $M$.

\item[c)] The planar representation has the schematic structure shown in Fig. \ref{fig:ccs}: the boundary components, represented by solid lines, lie either along the bottom or the top edge of the planar representation. Those along the bottom edge are called ``in'' boundary components, those along the top edge are called ``out'' boundary components. When there are no ``in/out'' boundary components, the bottom/top edge contains a single 0-cell. The dotted lines on the left and right, and the dotted lines between boundary components along the bottom and top edge, each represent one or more cut 1-cells, separated by 0-cells when there are more than one of them.

\item[d)] The simply connected planar region is made up of one or more 2-cells, separated by 1-cells and 0-cells when there are more than one of them.

\end{itemize}

\label{def:ccs}

\end{Definition}

\begin{Remark}
We will refer to the 0-cells and 1-cells that do not belong to a boundary component as internal or non-boundary 0-cells and 1-cells.
\end{Remark}

\begin{figure}[htbp] 
\centerline{\relabelbox 
\epsfysize 2.3cm
\epsfig{file=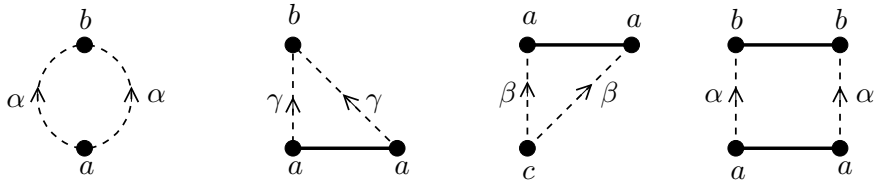,height=2.3cm}
\relabel{a}{$a$}
\relabel{b}{$b$}
\relabel{c}{$a$}
\adjustrelabel <-1pt,0pt> {d}{$a$}
\relabel{e}{$b$}
\relabel{f}{$c$}
\relabel{g}{$a$}
\adjustrelabel <0pt,1pt> {h}{$a$}
\relabel{i}{$a$}
\adjustrelabel <-1pt,0pt> {j}{$a$}
\relabel{k}{$b$}
\relabel{l}{$b$}
\relabel{X}{$\alpha$}
\relabel{Y}{$\alpha$}
\relabel{Z}{$\gamma$}
\relabel{W}{$\gamma$}
\adjustrelabel <-1pt,0pt> {T}{$\beta$}
\relabel{U}{$\beta$}
\relabel{R}{$\alpha$}
\relabel{S}{$\alpha$}
\endrelabelbox}
\caption{\label{fig:SDC} Examples of CCS's }
\end{figure}

To fix ideas we give some examples of cut cellular surfaces (Figure \ref{fig:SDC}), representing the sphere $S$, the disk $D$ (two versions with the boundary being ``in'' or ``out'' ) and the cylinder $C$. See \cite{bp1} for further examples and discussion.

{\bf Moves on CCS's.} By analogy with the Pachner moves on triangulated manifolds, we introduce moves for passing between different planar representations of the same surface. There are two types of move. 

\vskip 0.3cm

\begin{figure}[htbp] 
\centerline{\relabelbox 
\epsfxsize 8cm
\epsfig{file=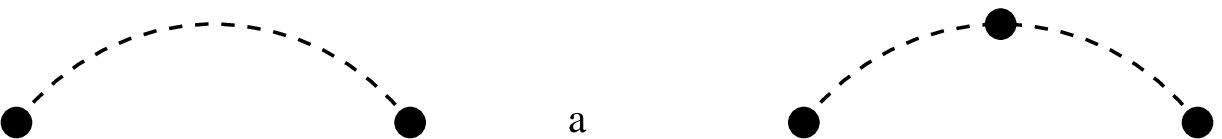,width=8cm}
\relabel{a}{$\longleftrightarrow$}
\endrelabelbox}
\caption{\label{fig:move1} Move I}
\end{figure}

{\bf Move I:} Introducing a 0-cell into a non-boundary 1-cell, thereby dividing it into two 1-cells, or conversely removing a 0-cell separating two 1-cells, to combine them into a single 1-cell (Figure \ref{fig:move1}). When this move is applied to a cut 1-cell, the 0-cell is introduced into or removed from both copies of the cut 1-cell in the planar representation. 

\begin{figure}[htbp] 
\centerline{\relabelbox 
\epsfxsize 13cm
\epsfig{file=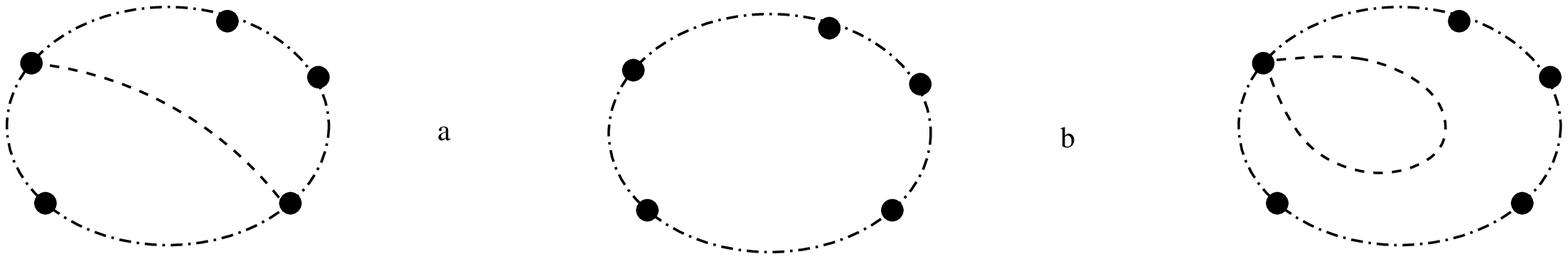,width=13cm}
\relabel{a}{$\longleftrightarrow$}
\relabel{b}{$\longleftrightarrow$}
\endrelabelbox}
\caption{\label{fig:move2} Move II}
\end{figure}

{\bf Move II:} Introducing a 1-cell into a 2-cell, thereby dividing it into two 2-cells, or conversely removing a 1-cell separating two 2-cells, to combine them into a single 2-cell (Figure \ref{fig:move2}).
In this figure we have used lines with dots and dashes for the 1-cells bounding the 2-cell to indicate that these are either boundary or internal 1-cells in the planar representation.

In \cite{bp1} it was shown that these moves generate the Pachner moves when $M$ is a triangulated surface without boundary.

\section{Invariants for CCS's from finite 2-group colourings}

\label{sec:inv}

We will be considering colourings of CCS's with finite crossed modules.

\begin{Definition}

A finite crossed module, or finite 2-group, $\G = (G,H, \partial, \triangleright)$ is given by:
\begin{itemize}
\item two finite groups $G$ and $H$
\item a group homomorphism $\partial : H \rightarrow G$
\item a left action $\triangleright$ of $G$ on $H$ by automorphisms
\end{itemize} 
such that, for all $h, h_1,h_2 \in H$ and $g \in G$:
\begin{align}
\partial(g \triangleright h) &= g \, \partial(h) \, g^{-1}\\
\partial(h_1) \triangleright h_2 &= h_1 h_2 h_1^{-1}
\end{align}
\label{def:2group}
\end{Definition}

\begin{Remark}An obvious class of examples is given by taking $G$ and $H$ to be the same, with  $\partial$ the identity, and $\tr$ given by conjugation. 
For any crossed module $\ker \partial$ is contained in the centre of $H$ and hence is abelian, since for $h\in \ker \partial$:  $ hfh^{-1}= \partial(h) \tr f = 1\tr f = f$.
A further class of examples comes from central extensions. Given a central extension of groups: 
$$
1\rightarrow K \rightarrow H \stackrel{\partial}{\rightarrow} G \rightarrow 1 
$$
one obtains a crossed module:
 $$
H \stackrel{\partial}{\longrightarrow} G
$$ 
with lifted action 
$$
g\tr h= fhf^{-1}
$$
where $f\in H$ is any element  such that $\partial(f)=g$. This action is well-defined 
because $K=\ker \partial$ is central in $H$.
\end{Remark}

Fix a finite crossed module ${\cal G} = (G,H, \partial, \triangleright)$. Given a CCS, $M$, we fix orientations on the 1-cells of $M$, specified as follows with respect to the planar representation:

\begin{itemize}

\item the boundary 1-cells are oriented from left to right
\item the cut 1-cells are oriented as chosen in Definition \ref{def:ccs} b)
\item the remaining internal 1-cells are oriented arbitrarily.

\end{itemize}

We also fix a basepoint (0-cell) in the boundary of each 2-cell.

\begin{Definition} A $\cal G$-colouring of $M$ is an assignment of an element $g_i\in G$ to each 1-cell labeled $i$ and of an element $h_A \in H$ to each 2-cell labeled $A$, such that, for each 2-cell in the planar representation, the following condition holds (which we call ``fake flatness'', in line with terminology from higher gauge theory in physics):

\begin{itemize}

\item if the 1-cells of the boundary of the 2-cell labeled $A$ are labeled $i_1, \dots i_k$, ordered in the anticlockwise direction starting at the basepoint, then
\begin{equation}
\prod_{j=1}^k g_{i_j}^{(-1)} = \partial(h_A)
\end{equation} 

where the factor is $g_{i_j}$ or $g_{i_j}^{-1}$, depending on whether or not the 1-cell $i_j$ is oriented compatibly with the positive orientation of the 2-cell.

\end{itemize}

\end{Definition}

\begin{figure}[htbp] 
\centerline{\relabelbox 
\epsfxsize 3cm
\epsfig{file=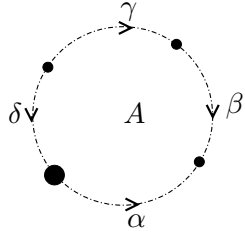,width=3cm}
\relabel{a}{$\alpha$}
\relabel{b}{$\beta$}
\relabel{c}{$\gamma$}
\relabel{d}{$\delta$}
\extralabel <-1.6cm,1.4cm> {$A$}
\endrelabelbox}
\caption{The fake flatness condition here is $g_\alpha g_\beta^{-1}g_\gamma^{-1}g_\delta =\partial (h_A)$.}
\label{fig:connection}
\end{figure}

See Figure \ref{fig:connection} for an example of the fake flatness condition. We have again used dots-and-dashes lines for the 1-cells to indicate that they can be either boundary or non-boundary 1-cells. The basepoint has been shown enlarged in the figure.

We can define invariants of CCS's, using $\cal G$-colourings. Choose elements of $G$, $g_1,\dots , g_n$, for the colouring of the ``in'' boundary components, and 
$g'_1,\dots , g'_m$, for the colouring of the ``out'' boundary components, ordering the boundary components from left to right in the planar representation. 
Let  $|G|$ denote the number of elements of the finite group $G$, $e$ denote the number of internal edges, i.e. 1-cells, and $v$ denote the number of internal vertices, i.e. 0-cells, of $M$. 
Let ${\rm Col}(g_1,\dots , g_n; g'_1,\dots , g'_m)$ denote the set of all $\cal G$-colourings of $M$ which have the given assignments on the boundary components.

\begin{Definition} The following invariants are defined for any choice of boundary colourings:
\begin{equation}
\left \langle g'_1,\dots , g'_m \left | Z_M \right |g_1,\dots , g_n   \right \rangle = \frac{|H|^{v-e}}{|G|^{\frac{m+n}{2}+v}} \, \# {\rm Col}(g_1,\dots , g_n; g'_1,\dots , g'_m).
\label{def:Z2grp}
\end{equation}
If $M$ has no ``in'' or ``out'' components we write the invariants as  $\left \langle {} \dots {} \left | Z_M \right |\emptyset  \right \rangle$ or $\left \langle \emptyset \left | Z_M \right | \dots    \right \rangle$. 
\label{def:2grpinv}
\end{Definition}

\begin{Remark}
When $M$ is a triangulated surface without boundary, these are the Yetter invariants \cite{yetter}. See \cite{FMPo} for an in-depth discussion of Yetter invariants.
\end{Remark}

\vskip 0.2cm
We now discuss in what sense these are invariants. First of all, we have:

\begin{Proposition}
The invariants $\left \langle g'_1,\dots , g'_m \left | Z_M \right |g_1,\dots , g_n   \right \rangle$ are unchanged under changes of orientation of the  internal 1-cells.
\label{prop:orientn}
\end{Proposition}
\begin{Proof} The number of internal vertices and edges is unchanged, and there is a bijection between the respective sets of colourings, given by replacing the element $g$ assigned to any  internal 1-cell by $g^{-1}$, when its orientation is reversed, thus guaranteeing that $h$ can be kept the same to satisfy the fake flatness condition.
\end{Proof}
\vskip 0.2cm

Likewise the choice of basepoints does not affect the invariant.

\begin{Proposition}
The invariants $\left \langle g'_1,\dots , g'_m \left | Z_M \right |g_1,\dots , g_n   \right \rangle$ are unchanged under a change of basepoint in any 2-cell.
\label{prop:start}
\end{Proposition}

\begin{Proof} The number of internal vertices and edges is unchanged, and there is a bijection between the respective sets of colourings, which only differ in the $H$-colouring of the 2-cell in question. The colourings $h$ and $h'$, corresponding to the first and second choice of basepoint respectively, are related by $h'= g^{-1} \triangleright h $ or equivalently $h= g \triangleright h' $, where $g$ is the ordered multiplication of the group colourings of the edges (taking into acccount orientation) that link the first basepoint to the second basepoint going round in the anticlockwise direction. This indeed establishes a bijection between the two sets of colourings, since the action of $G$ on $H$ is by automorphisms. The fake flatness condition for the first basepoint may be written as $\partial(h)=gk$, where $k$ represents the ordered multiplication of the group colourings of the edges (taking into acccount orientation) that link the second basepoint to the first basepoint going round in the anticlockwise direction. The fake flatness condition for the second basepoint then follows: $
\partial(h') = \partial(g^{-1} \triangleright h) = g^{-1} \partial(h) g = kg $.
\end{Proof}

\vskip 0.1cm

\begin{figure}[htbp] 
\centerline{\relabelbox 
\epsfxsize 11.5cm
\epsfig{file=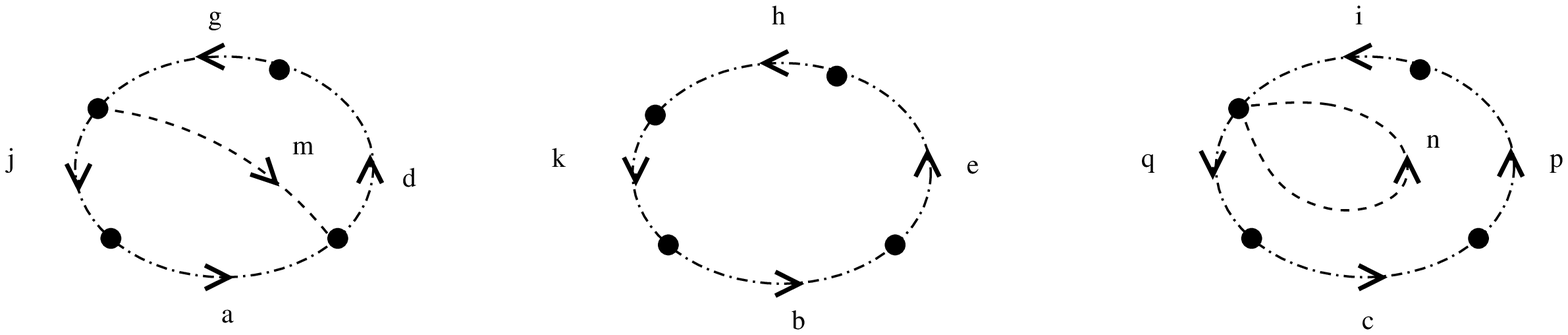,width=13.5cm}
\relabel{a}{$k_1$}
\relabel{b}{$k_1$}
\relabel{c}{$k_1$}
\relabel{d}{$k_2$}
\relabel{e}{$k_2$}
\relabel{p}{$k_2$}
\relabel{g}{$k_3$}
\relabel{h}{$k_3$}
\relabel{i}{$k_3$}
\relabel{j}{$k_4$}
\relabel{k}{$k_4$}
\relabel{q}{$k_4$}
\adjustrelabel <1pt,-8pt> {m}{$g$}
\adjustrelabel <-3pt,-1pt> {n}{$g$}
\extralabel <-12.2cm,1cm> {$h_1$}
\extralabel <-11.5cm,1.8cm> {$h_2$}
\extralabel <-7cm,1.2cm> {$h$}
\extralabel <-2.4cm,1.4cm> {$h_2$}
\extralabel <-1.8cm,.7cm> {$h_1$}
\endrelabelbox}
\caption{Colourings of $M$ (in the middle) and $M'$ for Move II}
\label{fig:pfmv2}
\end{figure}

More importantly we have:

\begin{Theorem}
The invariants $\left \langle g'_1,\dots , g'_m \left | Z_M \right |g_1,\dots , g_n   \right \rangle$ are unchanged under moves I and II.
\end{Theorem}
\begin{Proof} Suppose $M$ and $M'$ are related by a move I. Fix a $\cal G$-colouring for $M$ that assigns $g$ to the 1-cell displayed on the left in Figure \ref{fig:pfmv1}. Keeping the assignments of all other 1-cells and 2-cells the same, for $M'$ on the right there are $|G|$ compatible $\cal G$-colourings, since we can choose one assignment, e.g. $j$, freely in $G$ and the other assignment $k$ is then determined (for the orientations as shown in Figure \ref{fig:pfmv1}, we have $g=jk$, i.e. $k=j^{-1}g$). Since $M'$ has both an extra internal vertex and an extra internal edge compared to $M$, the exponent of $|H|$ in (\ref{def:Z2grp}) is unchanged, and the increase in the exponent of $|G|$ in the denominator is cancelled by the factor $|G|$ relating the respective number of colourings. Thus the invariants (\ref{def:Z2grp}) are the same for $M$ and $M'$.

\begin{figure}[htbp] 
\centerline{\relabelbox 
\epsfxsize 8cm
\epsfig{file=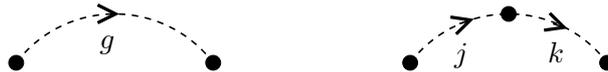,width=8cm}
\relabel{g}{$g$}
\relabel{j}{$j$}
\relabel{k}{$k$}
\endrelabelbox}
\caption{Colourings of $M$ and $M'$ for Move I}
\label{fig:pfmv1}
\end{figure}

Suppose $M$ and $M'$ are related by a move II, where $M'$ has an extra internal 1-cell compared to $M$, dividing a 2-cell in $M$ into two 2-cells in $M'$. Using basepoint invariance we may choose the basepoints of the two 2-cells in $M'$ to coincide, and we may choose this same 0-cell as the basepoint of the 2-cell in $M$. Using invariance under change of orientation, we may choose the extra 1-cell in $M'$ to be oriented so as to have the basepoint as its starting point  (see Figure \ref{fig:pfmv2}, where the basepoint for all 2-cells is the starting point of the 1-cell labelled $k_4$).

Fix a $\G$-colouring of $M$ that assigns $h\in H$ to the 2-cell we are considering. Keeping the assignments of all other 1- and 2-cells the same, there are $|H|$ corresponding $\G$-colourings of $M'$, since we may choose freely an element $h_2\in H$ to assign to the 2-cell, say on the left as we follow the subdividing 1-cell in the direction of its orientation, which then determines uniquely the assignment of 
$h_1=h h_2^{-1}$ to the other 2-cell, and the assignment of an element $g\in G$ to the subdividing 1-cell, by using the fake flatness condition in either 2-cell. These assignments are compatible with fake flatness, since imposing fake flatness implies $\partial(h_1)=\partial(hh_2^{-1})= \partial(h)\partial(h_2^{-1})$, which is necessary. Indeed, taking $M'$ on the left in Figure \ref{fig:pfmv2} as an example, $\partial(h_1)=k_4k_1g^{-1}$ and $\partial(h) \partial(h_2^{-1}) = k_4k_1k_2k_3 \, . \, k_3^{-1}k_2^{-1}g^{-1}$ are the same.

Conversely, given a $\G$-colouring of $M'$ which assigns $h_1$ and $h_2$ to the left and right 2-cell respectively, there is a compatible $\G$-colouring of $M$ which assigns $h=h_1h_2$ to the undivided 2-cell and agrees with the $\G$-colouring of $M'$ elsewhere. There are $|H|$ possible $\G$-colourings of $M'$ which give the same $h\in H$, namely $h_1'=h_1h'$ and $h_2'=(h')^{-1}h_2$ for any $h'\in H$.

$M'$ has the same number of internal vertices as $M$ and one extra internal 1-cell. Thus the exponent of $|G|$ in (\ref{def:Z2grp}) is the same for both $M$ and $M'$, and the increase in the number of colourings for $M'$ by a factor $|H|$ is cancelled by the extra factor $|H|^{-1}$ in (\ref{def:Z2grp}) coming from the extra internal 1-cell. Thus the invariants are the same for $M$ and $M'$.
\end{Proof}

\vskip 0.3cm

\section{Examples}
\label{sec:exp}

\begin{figure}[htbp] 
\centerline{\relabelbox 
\epsfysize 2.8cm
\epsfig{file=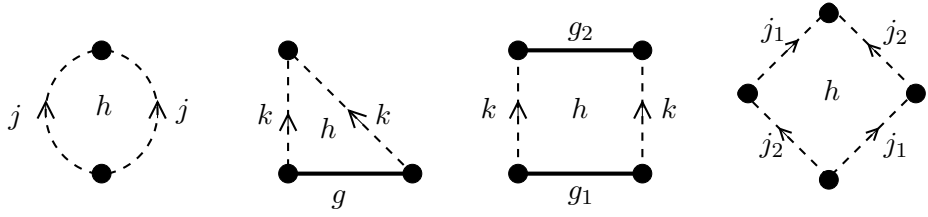,height=2.8cm}
\relabel{D}{$j_1$}
\adjustrelabel <-3pt,1pt> {B}{$j_1$}
\adjustrelabel <-2pt,0pt> {A}{$j_2$}
\relabel{C}{$j_2$}
\relabel{X}{$j$}
\relabel{Y}{$j$}
\relabel{Z}{$k$}
\relabel{W}{$k$}
\relabel{R}{$k$}
\relabel{S}{$k$}
\adjustrelabel <0pt,3pt> {T}{$g$}
\adjustrelabel <0pt,2pt> {U}{$g_1$}
\adjustrelabel <0pt,1pt> {h}{$g_2$}
\adjustrelabel <0pt,-1pt> {a}{$h$}
\adjustrelabel <-2pt,-3pt> {b}{$h$}
\adjustrelabel <-1pt,0pt> {c}{$h$}
\adjustrelabel <-2pt,0pt> {d}{$h$}
\endrelabelbox}
\caption{\label{fig:inv-exp1} $\cal G$-colourings for the sphere, disk, cylinder and torus. }
\end{figure}

In this section we calculate the invariant for some simple examples.  Let $K$ and $A$ denote the kernel and image of $\partial$, respectively. In Figure \ref{fig:inv-exp1} below we choose the basepoint to be the bottom 0-cell and on the left, if there is a choice. 

Starting with the disk $D$, the fake flatness condition is $\partial(h) = gkk^{-1}=g$. There are two ways to calculate the overall number of colourings, allowing arbitrary $g\in G$. The first is to fix the 2-cell coloring $h$. Then $g$ is determined through the condition and the number of colourings is $|H|$. The second is to fix the 1-cell coloring $g$. If $g \in A$, then one has $|K|$ possible values for $h$ and the number of colourings is $|A|\, |K|$, which is compatible with the previous result, since $|H|=|A|\,|K|$ from group theory. If $g \notin A$, no colourings are possible. The disk has one internal vertex ($v=1$), one internal edge ($e=1$), and one boundary component ($m=0,\,n=1$). Thus we have:
$$
\langle \emptyset |Z_D| g \rangle = \frac{1}{|G|^{1/2}} |K| {\cal D}(g)\quad {\rm where} \quad {\cal D}(g) :=\left\{ \begin{array}{cc} 1, & g \in A \\ 0, & g\notin A \end{array} \right . 
$$

For the sphere $S$, the fake flatness condition is $\partial(h)=jj^{-1}=1$, i.e. we have $h\in K$, and $j$ is arbitrary in $G$. Hence the number of colourings is 
$|K|\, |G|$, which together with $v=2,\, e=1,\, m=n=0$, leads to:
$$
\langle \emptyset |Z_S| \emptyset \rangle = \frac{|H|}{|G|^2} |K| |G|= \frac{|H||K|}{|G|} \,
$$

For the cylinder $C$, the fake flatness condition gives $\partial(h) = g_1 k g_2^{-1} k^{-1}$. Note that, if $A$ is the trivial group with one element, this condition expresses that $g_1$ and $g_2$ are conjugate to each other, since it is equivalent to: $g_2= k^{-1} g_1 k$. We will have more to say about the relation between $g_1$ and $g_2$ in section \ref{sec:glue}. For $C$ we have $v=0,\, e=1, \, m=n=1$, and hence
$$
\langle g_2 |Z_C| g_1 \rangle = \frac{1}{|H| |G|} {\cal C}(g_1,g_2) 
$$
where
\begin{equation}
{\cal C}(g_1,g_2)=\# \{ (h,k) \in H \times G : \partial(h) = g_1 k g_2^{-1} k^{-1} \}
\label{eq:Cdef}
\end{equation}

Finally, for the torus $T$, the fake flatness condition is $\partial(h)=j_1 j_2 j_1^{-1} j_2^{-1} $, meaning that $h$ has to lie in the preimage under $\partial$ of the commutator subgroup of $G$. Since $v=1,\, e= 2,\, m=n=0$, we have:
\begin{align}
\langle \emptyset |Z_T| \emptyset \rangle &= \frac{\#\{(h,g_1,g_2) \in H \times G^2 : \partial(h)=j_1 j_2 j_1^{-1} j_2^{-1}\}}{|G| |H|} 
\label{eq:Ztorus}
\end{align}

\section{Gluing formula and TQFT}
\label{sec:glue}

\begin{figure}[htbp] 
\centerline{\relabelbox 
\epsfxsize 10cm
\epsfig{file=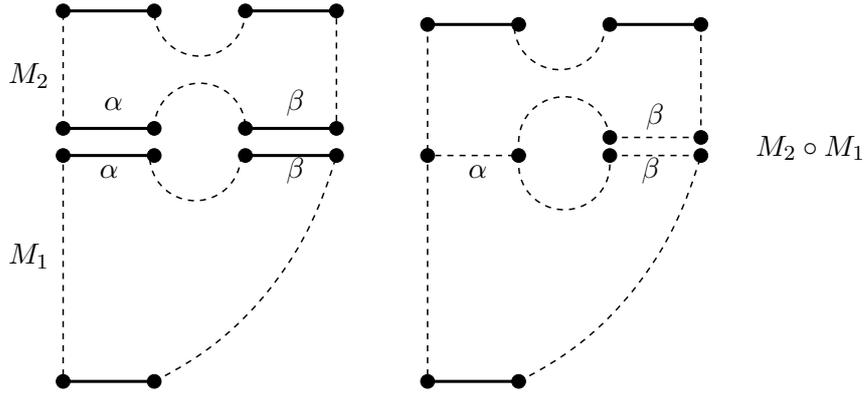,width=10cm}
\relabel{a}{$M_1$}
\relabel{b}{$M_2$}
\relabel{c}{$M_2\circ M_1$}
\relabel{h}{$\alpha$}
\relabel{i}{$\alpha$}
\relabel{e}{$\alpha$}
\relabel{j}{$\beta$}
\relabel{d}{$\beta$}
\relabel{f}{$\beta$}
\relabel{g}{$\beta$}
\endrelabelbox}
\caption{Gluing or composition of two CCS's}
\label{fig:gluing}
\end{figure}

Following our approach in \cite{bp1}, in order to glue two surfaces $M_1$ and $M_2$ with matching boundaries, we adopt the following procedure:
we identify the shared boundary component furthest to the left (labeled $\alpha$ in Figure \ref{fig:gluing}), and the remaining shared boundary components (just one in the Figure, labeled $\beta$) become cut 1-cells in the boundary of the planar representation of a new CCS that we denote by $M_2 \circ M_1$.

Suppose we have $M_1$ with $n$ incoming boundary components and $m>0$ outgoing boundary components, and $M_2$ with $m$ incoming boundary components and $p$ outgoing boundary components. Fixing the colourings of the ``in'' boundary components of $M_1$ and the ``out'' boundary components of $M_2$ the colourings of $M_2\circ M_1$ allow a priori any choice for the colourings of the $m$ intermediate 1-cell components. Thus we arrive at the following property for the invariants.

\begin{Proposition} (Gluing formula) For any $g_1, \dots , g_n, i_1,\dots ,i_p\in G$, we have:
\begin{eqnarray}
\lefteqn{\left \langle i_1,\dots , i_p \left | Z_{M_2\circ M_1} \right |g_1,\dots , g_n   \right \rangle =} \nonumber \\  
& &  \sum_{j_1,\dots , j_m\in G }  
\left \langle i_1,\dots , i_p \left | Z_{M_2} \right | j_1,\dots , j_m  \right \rangle  
\left \langle j_1,\dots , j_m \left | Z_{M_1} \right |g_1,\dots , g_n   \right \rangle 
\label{eq:gluing}
\end{eqnarray}
\end{Proposition}
\begin{Proof}
Since the number of colourings match on both sides of the equation, it remains to check the other factors. Each of the $2m$ boundary components that are glued in $M_1$ and $M_2$  
gives rise to a factor $\frac{1}{|G|^{1/2}}$ in (\ref{def:Z2grp}). After gluing $M_2\circ M_1$ has instead $m$ extra internal vertices, each of which gives a factor $\frac{1}{|G|}$. The factors of $|H|$ are also the same, since for  $M_2\circ M_1$ the $m$ additional internal vertices and the $m$ additional internal edges cancel in the exponent of $|H|$ in (\ref{def:Z2grp}).
\end{Proof} 
\vskip 0.2cm

The gluing formula enables us to construct a natural TQFT - see our previous article \cite{bp1} for an introduction to the notion of TQFT. We assign to each incoming or outgoing boundary of a CCS $M$, a vector space $V_{\rm in}$ or $V_{\rm out}$ over $\mathbb{R}$, whose basis consists of all $G$-colourings of the boundary components \cite{bp1}. The basis elements are written $|\,g_1,\dots , g_n  \left. \right \rangle$ or  
$\left \langle \right. i_1,\dots , i_m\,  |$, and the dimension of $V_{\rm in}$ and $V_{\rm out}$ is $|G|^n$ and $|G|^m$ respectively. To the CCS itself we assign the linear transformation $Z_M$ from $V_{\rm in}$ to $V_{\rm out}$, whose matrix elements with respect to these two bases are given by:

$$
 \left \langle i_1,\dots , i_m \left | Z_M \right |g_1,\dots , g_n   \right \rangle
$$

Thus from the gluing formula (\ref{eq:gluing}) we have the following fundamental result:
\begin{Proposition} (TQFT property) For any $M_1$ and $M_2$ such that $M_2\circ M_1$ is defined, we have:
\begin{equation}
Z_{M_2\circ M_1} = Z_{M_2} \circ Z_{M_1}.
\label{tqft}
\end{equation}
\end{Proposition}

There is an important corollary of (\ref{tqft}), which expresses that the cylinder $C$ has assigned to it an idempotent (since $C\circ C$ and $C$ are related by moves I and II, we have 
$Z_{C\circ C} =Z_C$):
\begin{Corollary} For the cylinder $C$, $Z_C$ satisfies
$$
Z_C \circ Z_C = Z_C 
$$
\end{Corollary}

In terms of the function ${\cal C}$ defined in \eqref{eq:Cdef}, this result is equivalent to
\begin{equation}
\sum_{i \in G} {\cal C} \left(g,i\right) {\cal C} \left(i,j\right)  = |H| |G| \cdot {\cal C} \left(g,j\right)\, ,
\label{eq:cyl}
\end{equation}
for any $g,j\in G$.

We also give an algebraic proof of (\ref{eq:cyl}), which will be useful in what follows. First we define an equivalence relation in $G$.
\begin{Definition} 
We say that two elements $g_1$ and $g_2$ of $G$ are 2-conjugate in $\G$, denoted $g_1 \sim g_2 $, iff 
$$
{\cal C}(g_1,g_2)\neq 0 \, .
$$.
\end{Definition}

\begin{Proposition} 
2-conjugacy is an equivalence relation.
\end{Proposition}

\begin{Proof}
Let $W(g_1,g_2)$ denote the set $\{ (h,k) \in H \times G : \partial(h) = g_1 k g_2^{-1} k^{-1} \}$, which has cardinality ${\cal C}(g_1,g_2)$. Then $g_1 \sim g_2 $ iff $W(g_1,g_2)\neq \emptyset$.

$\sim$ is reflexive, since $(1_H,1_G)\in W(g,g)$ for every $g\in G$.

$\sim$ is symmetric: if $(h,k) \in W(g_1,g_2)$, then $(k^{-1}\tr h^{-1},k^{-1}) \in W(g_2,g_1)$, since 
\begin{eqnarray*}
\partial(k^{-1}\tr h^{-1}) & = & k^{-1} \partial(h^{-1})k \\
&=&  k^{-1} (kg_2k^{-1}g_1^{-1})k \\
&=& g_2 k^{-1} g_1^{-1} k
\end{eqnarray*}

$\sim$ is transitive: if $(h,k) \in W(g_1,g_2)$  and $(h',k') \in W(g_2,g_3)$, then we have
$(h(k\tr h'),kk') \in W(g_1,g_3)$, since
\begin{eqnarray*}
\partial(h(k\tr h')) & = & \partial(h) \partial(k\tr h') \\
&=& \partial(h) k \partial( h') k^{-1} \\
& = & g_1 k g_2^{-1} k^{-1} k  (g_2 k' g_3^{-1} k'^{-1}) k^{-1} \\
&=& g_1 (kk') g_3^{-1} (kk')^{-1}
\end{eqnarray*}
\end{Proof}
\vskip 0.2cm

Using the sets $W(g_1,g_2)$ introduced in the the previous proof, we can also show the following symmetry.
\begin{Proposition} 
For all $g_1,\, g_2\in G$, we have ${\cal C}(g_1,g_2)={\cal C}(g_2,g_1)$.
\end{Proposition}
\begin{Proof}
We establish a bijection $W(g_1,g_2) \underset{\alpha}{\overset{\beta}{\leftrightarrows}}  W(g_2,g_1)$ by defining
$$
\alpha(h,k)=(k^{-1}\tr h^{-1}, k^{-1}) \qquad \beta(h',k') = (k'^{-1}\tr h'^{-1},k'^{-1})
$$
From the previous proof, $\alpha$ is well-defined, i.e. 
$\partial(k^{-1}\tr h^{-1})  g_2(k^{-1})g_1^{-1} k$,
and likewise $\beta$ is well-defined. 
$\alpha$ and $\beta$ constitute a bijection since 
\begin{eqnarray*}
(\beta\circ\alpha)(h,k) & = & \beta(k^{-1}\tr h^{-1}, k^{-1}) \\
& = & (k\tr (k^{-1}\tr h),k) \\
& = & (h,k)
\end{eqnarray*}
and likewise $(\alpha\circ\beta)(h',k')=(h',k')$. Thus $W(g_1,g_2)$ and  $W(g_2,g_1)$ are isomorphic, and hence their cardinality is the same. 
\end{Proof}

Using analogous methods we have the following result.
\begin{Proposition} 
If $g_1 \sim g_2$, then ${\cal C}(g_1,g_2)={\cal C}(g_1,g_1)$.
\end{Proposition}
\begin{Proof}
Since $g_1 \sim g_2$, there exists a pair $(h,k)\in H\times G$ such that 
$\partial(h) = g_1 k g_2^{-1} k^{-1} $.

We establish a bijection $W(g_1,g_2) \underset{\alpha}{\overset{\beta}{\leftrightarrows}}  
W(g_1,g_1)$ by defining
$$
\alpha(h',k')=(h'(k'k^{-1})\tr h^{-1}, k'k^{-1}) \qquad 
\beta(h'',k'') = (h''(k''\tr h),k''k)
$$
$\alpha$ is well-defined, since 
\begin{eqnarray*}
\partial(k^{-1}\tr h^{-1}) & = & k^{-1} \partial(h^{-1}) k  \\
& = & k^{-1} kg_2k^{-1}g_1^{-1}k \\
& = & g_2(k^{-1})g_1^{-1} (k^{-1})^{-1}
\end{eqnarray*}
and likewise $\beta$ is well-defined. 
$\alpha$ and $\beta$ constitute a bijection since 
\begin{eqnarray*}
(\beta\circ\alpha)(h,k) & = & \beta(k^{-1}\tr h^{-1}, k^{-1}) \\
& = & (k\tr (k^{-1}\tr h),k) \\
& = & (h,k)
\end{eqnarray*}
and likewise $(\alpha\circ\beta)(h',k')=(h',k')$. Thus $W(g_1,g_2)$ and  $W(g_1,g_1)$ are isomorphic, and hence their cardinality is the same. 
\end{Proof}
\vskip 0.2cm

Using the previous proposition for the non-zero terms on the l.h.s. of (\ref{eq:cyl}), we have $i\sim j$, hence ${\cal C}(i,j)={\cal C}(j,j)$. Also
$g\sim i\sim j$, and therefore ${\cal C}(g,j)={\cal C}(j,j)$. Thus (\ref{eq:cyl}) is equivalent to
$$
\sum_{i\in G}{\cal C}(g,i)=|H||G|
$$
which clearly holds since, fixing $g$, every pair $(h,k)\in H 	\times G$ belongs to precisely one set of the form $W(g,i)$ with $g$ fixed.

\begin{Remark}
If we denote the 2-conjugacy class of $g \in G$ by $\bar g$, we get an equation for the number of elements of $\bar g$:
\begin{equation}
\# \bar{g} = \frac{|G| |H|}{{\cal C}(g,g)} \,,
\end{equation}
since ${\cal C}(g,g_1) = {\cal C}(g,g)$ for all $g_1 \in \bar{g}$, and the number of non-zero terms in the sum on the l.h.s of Eq.~\eqref{eq:cyl} is $\# \bar{g}$.
Let ${\rm 2ConjClass}(G)$ denote the set of 2-conjugacy classes of $G$. Then its cardinality is given by the following equation
\begin{equation}
\label{eq:equiv}
\#{\rm 2ConjClass}(G) = \frac{1}{|G|^2 |H|^2} \sum_{g,g_1 \in G} {\cal C}(g,g_1)^2\,.
\end{equation} 
This is clear since the double sum on the r.h.s. decomposes into double sums where $g,g_1$ both belong to the same 2-conjugacy class $\bar g$. Restricting to these terms for a specific 2-conjugacy class 
$\bar g$, the r.h.s. of Eq.~\eqref{eq:equiv} becomes:
\begin{equation}
\frac{1}{|G|^2 |H|^2}  \sum_{g,g_1 \in \bar g} {\cal C}(g,g_1)^2 = \frac{(\# \bar{g})^2 {\cal C}(g,g)^2}{|G|^2 |H|^2} = 1 \,,
\end{equation}
and collecting the contributions from each class, we obtain equation \eqref{eq:equiv}.

\end{Remark}

Our final result is analogous to a proposition obtained in \cite{bp1}. Consider the invariant for the torus (\ref{eq:Ztorus}), which may be rewritten, using \eqref{eq:cyl}, as follows: 
\begin{align}
\langle \emptyset |Z_T| \emptyset \rangle &= \frac{\#\{(h,g_1,g_2) \in H \times G^2 : \partial(h)=g_1 g_2 g_1^{-1} g_2^{-1}\}}{|G| |H|} \nonumber \\
&= \frac{\sum_{g_1 \in G} {\cal C}(g_1,g_1)}{|G| |H|} \nonumber\\
&= \frac{\sum_{g_1,g_2 \in G} {\cal C}(g_1,g_2)^2}{|G|^2 |H|^2} \, ,
\label{eq:T=2C}
\end{align}
This equation reflects the topological fact that the torus is obtained by gluing two cylinders together. See \cite{bp1} where this point was developed. 

For a finite group (not a 2-group) $G$, its commuting fraction is defined to be the ratio 
$$
\frac{\#\{(g_1,g_2) \in G^2 : g_1 g_2 = g_2 g_1 \}}{|G|^2} \, ,
$$
i.e. the ratio of the number of commuting pairs of elements over the number of all pairs of elements.
Here we define an analogous fraction for a 2-group $\G$, namely the \emph{generalized commuting fraction}  of $\G$
\begin{equation}
\label{eq:commutfrac}
\frac{\#\{(h,g_1,g_2) \in H \times G^2 : \partial(h)=g_1 g_2 g_1^{-1} g_2^{-1}\}}{|H| |G|^2} \,.
\end{equation}
By combining  equations \eqref{eq:equiv} and \eqref{eq:T=2C} with definition (\ref{eq:commutfrac}), we obtain the following generalization of proposition 6.5 in \cite{bp1}.
\begin{Proposition}
The number of 2-conjugacy classes of $G$ in $\G$ is equal to the generalized commuting fraction of $\cal G$ times the order of $G$.
\label{prop:gcf}
\end{Proposition}

\section{Conclusions and Final Remarks}
\label{sec:conc}

Viewing our results from the 2-group theory perspective, we have been led to introduce a function $\cal C$, taking values in the non-negative integers, which depends on two $G$-elements. The function $\cal C$ defines an equivalence relation, 2-conjugacy, between $G$-elements, namely $g_1$ and $g_2$ are 2-conjugate iff ${\cal C}(g_1,g_2)\neq 0$, but also gives a ``measure of the equivalence'' between the elements $g_1$ and $g_2$ by counting the number of pairs $(h,k) \in H \times G$ such that $\partial(h) = g_1 k g_2^{-1} k^{-1}$. We have derived properties of $\cal C$ by using topological reasoning, leading us to define the  generalized commuting fraction for a 2-group, which we proved to have a property analogous to a property of the standard commuting fraction of a finite group. It should be possible to obtain many further results in the theory of finite 2-groups using a similar topological approach. 

Using colourings of cut cellular surfaces with elements of a finite 2-group $\cal G$,  we have found not only invariants for these surfaces, but also a TQFT setting for the invariants. Interestingly these TQFT's do not naturally fit into the standard framework of Atiyah's axioms for TQFT (see \cite{at}  and the study by Abrams of 2-dimensional TQFT's \cite{ab}), since $Z_{ C}$ for the cylinder is an idempotent, not necessarily the identity. We note that the eigenvalue 1 eigenvectors of $Z_{C}$ are of the form $g_1 + g_2 + \dots + g_k$ where the sum runs over all elements of a 2-conjugacy class in $G$. 

As already alluded to, there is an interpretation of these invariants in terms of higher gauge theory based on a finite 2-group $\G$ which we now sketch. First we look at the invariant of our previous article \cite{bp1} from the point of view of ordinary gauge theory based on a finite group $G$. In this context we are interested in the moduli space of flat $G$-connections modulo gauge transformations, or rather the corresponding groupoid whose objects are flat $G$-connections and morphisms are gauge transformations. Flat $G$-connections correspond to flat $G$-colourings of the 1-cells of the surface $M$, and the gauge transformations are given by assignments of elements of $G$ to the 0-cells of $M$. For a surface without boundary $M$, the invariant of \cite{bp1} 
$$
\langle \emptyset |Z_M| \emptyset \rangle
= 
\frac{\#\, {\rm Flat}\,  G{\text -}{\rm colourings}}{|G|^v}
$$
where $v$ denotes the number of 0-cells, can be understood as the {\em groupoid cardinality}\footnote{ The groupoid cardinality can be thought of as counting the objects of a groupoid taking into account the number of isomorphisms each object has with other objects. For a nice introduction to the notion of groupoid cardinality, see \cite{tao}, and for the more general concept of the Euler characteristic of a category, see \cite{lei}. } of the groupoid of flat $G$-connections on $M$. In the case of a groupoid coming from the action of a finite group $\widetilde{G}$ on a finite set $S$, the groupoid cardinality is simply the quotient of the respective cardinalities: $|S|/|\widetilde{G}|$.

In higher gauge theory an analogous picture is emerging, in work by one of us with J. Morton \cite{mp1,mp2}. The higher connections are given by $\G$-colourings of the 1- and 2-cells of $M$ satisfying the fake flatness condition. There are two different types of gauge transformation between these connections, corresponding to assignments of $G$ elements to the 0-cells of $M$ as well as assignments of $H$ elements to the 1-cells of $M$ (taken to be without boundary). In addition there are higher-level transformations between gauge transformations given by assignments of $H$ elements to the 0-cells of $M$. A satisfying description of all this is in terms of a double groupoid, i.e. a higher algebraic structure having objects, two types of morphisms between objects called horizontal and vertical, and higher morphisms called squares between the morphisms, all morphisms being suitably invertible. The invariant (\ref{def:Z2grp}), written as follows:
 $$
\langle \emptyset |Z_M| \emptyset \rangle
= 
\frac{\#\, {\rm Fake}\, {\rm Flat}\,  \G{\text -}{\rm colourings}\, .	\, {|H|^v}}{|G|^v\, |H|^e}
$$
can thus naturally be viewed as the ``double groupoid cardinality'' of the double groupoid of higher connections. This perspective will be explored in more detail elsewhere.

\section{Acknowledgments}

This article is based on a research project carried out by Diogo Bragança under the supervision of  Roger Picken. The authors are grateful to the {\em Fundação Calouste Gulbenkian} for supporting this project through the programme {\em Novos Talentos em Matemática}, which aims to stimulate undergraduate research in mathematics.  
Roger Picken is grateful to Dan Christensen and Jeffrey Morton for useful discussions and suggestions. This work was funded in part by the Center for Mathematical Analysis, Geometry and Dynamical Systems (CAMGSD),
Instituto Superior T\'ecnico, Universidade de Lisboa, through the project UID/MAT/04459/2013 of the {\em Fundação para a Ciência e a Tecnologia} (FCT, Portugal).


\begin{thebibliography}{99}

\bibitem{ab} L. Abrams, ``Two-dimensional topological quantum field theories and Frobenius algebras'', 
\emph{J. Knot Theory Ramifications} Vol.~5, No.~5 (1996), pp. 569–587. 
\url{https://doi.org/10.1142/S0218216596000333}

\bibitem{at} M. Atiyah, ``Topological quantum field theories'', \emph{Inst. Hautes Études Sci. Publ. Math.} No. 68 (1988), pp. 175–186 (1989). 
\url{http://www.numdam.org/item?id=PMIHES_1988__68__175_0}

\bibitem{bp1} D. Bragança and R. Picken,``Invariants and TQFT’s for cut cellular surfaces from finite groups'',  \emph{Bol. Soc. Port. Mat.} Vol.~74 (2016), pp. 17-44.


\bibitem{dw} R. Dijkgraaf and E. Witten, ``Topological Gauge Theories and Group Cohomology'', \emph{Comm.  Math. Phys.} Vol.~129, No.~2 (1990), pp. 393–429. 
\url{https://projecteuclid.org/euclid.cmp/1104180750}

\bibitem{FMPo} J. Faria Martins and T. Porter, ``On Yetter's invariant and an extension of the Dijkgraaf-Witten invariant to categorical groups'', \emph{Theory Appl. Categ.} Vol.~18, No. 4 (2007), pp. 118–150. 
\url{http://www.tac.mta.ca/tac/volumes/18/4/18-04abs.html}

\bibitem{lei} T. Leinster, ``The Euler characteristic of a category'', \emph{ Doc. Math.} Vol. 13 (2008), pp. 21-49. 
\url{https://www.math.uni-bielefeld.de/documenta/vol-13/02.pdf}

\bibitem{mp1} J. C. Morton and R. Picken, ``Transformation double categories associated to 2-group actions'', \emph{Theory Appl. Categ.} Vol. 30, Paper No. 43 (2015), pp. 1429–1468. 
\url{http://www.tac.mta.ca/tac/volumes/30/43/30-43abs.html}

\bibitem{mp2} J. C. Morton and R. Picken, ``2-Group Actions on Moduli Spaces of Higher Gauge Theory'' (in preparation).

\bibitem{tao} T. Tao, ``Counting objects up to isomorphism: groupoid cardinality'', What's New (2017).\\
\href{https://terrytao.wordpress.com/2017/04/13/counting-objects-up-to-isomorphism-groupoid-cardinality/}{\nolinkurl{https://terrytao.wordpress.com/2017/04/13/counting-objects-}} \\
\href{https://terrytao.wordpress.com/2017/04/13/counting-objects-up-to-isomorphism-groupoid-cardinality/}{\nolinkurl{up-to-isomorphism-groupoid-cardinality/} }


\bibitem{yetter} D. N. Yetter, ``TQFT’s from Homotopy 2-Types'', \emph{J. Knot Theory Ramifications} Vol.~2, No.~1 (1993), pp. 113–123.
\url{https://doi.org/10.1142/S0218216593000076}



\end{thebibliography}
\end{document}